\documentclass[10pt,a4paper]{amsart}
\usepackage[english]{babel}
\usepackage{amssymb,xypic,amscd}
\CompileMatrices
\newtheorem{defn}{Definition}

\theoremstyle{remark}
\newtheorem{rem}[defn]{Remark}
\theoremstyle{remark}

\numberwithin{equation}{section} \numberwithin{defn}{section}

\newcommand\ed{\operatorname{End}}
\newcommand\tr{\operatorname{Tr}}

\begin{document}

\title[Linearity of Tate's Trace]{A negative answer to the
question \\ of the linearity of Tate's Trace \\ for the sum of two
endomorphisms}
\author{Julia Ramos Gonz\'alez (*)  \\ Fernando Pablos Romo (**)}

\address{ Departamento de Matem\'aticas, Universidad de Salamanca, Plaza de la Merced 1-4, 37008 Salamanca, Espa\~na}
 \email{ (*) julk@usal.es}
 \email{ (**) fpablos@usal.es}

\thanks{
\\ This work is partially supported by a collaboration-fellowship of the Spanish Government (*), and by the
DGESYC research contract no. MTM2009-11393 (**)}

\maketitle

\begin{abstract} The aim of this note is to solve a problem proposed by J. Tate in 1968 by offering a counter-example
of the linearity of the trace for the sum of two finite potent
operators on an infinite-dimensional vector space.
\end{abstract}

\bigskip

\setcounter{tocdepth}1

\tableofcontents
\bigskip

\section{Introduction}

    Let $k$ be a fixed ground field and $V$ a vector
space over $k$. If we consider an endomorphism $\varphi$ of $V$,
according to \cite{Ta} we say that $\varphi$ is ``finite-potent"
if $\varphi^n V$ is finite dimensional for some $n$, and a trace
$\tr_V(\varphi) \in k$ may be defined, having the following
properties:
\begin{enumerate}
\item If $V$ is finite dimensional, then $\tr_V(\varphi)$ is the
ordinary trace. \item If $W$ is a subspace of $V$ such that
$\varphi W \subset W$, then $$\tr_V(\varphi) = \tr_W(\varphi) +
\tr_{V/W}(\varphi)\, .$$ \item If $\varphi$ is nilpotent, then
$\tr_V(\varphi) = 0$. \item If $F$ is a ``finite-potent'' subspace
of $\ed (V)$ (i.e., if there exists an $n$ such that for any
family of $n$ elements $\varphi_1,\dots ,\varphi_n \in F$, the
space $\varphi_1\dots \varphi_n V$ is finite dimensional), then
$\tr_V\colon F\longrightarrow k$ is $k$-linear. \item If $f\colon
V'\to V$ and $g\colon V\to V'$ are $k$-linear and $f\circ g$ is
finite potent, then $g\circ f$ is finite potent, and $$\tr_V
(f\circ g) = \tr_{V'} (g\circ f)\, .$$
\end{enumerate}
\begin{rem}\label{r:chartraces} Properties (1), (2) and (3)
characterize traces, because if $W$ is a finite dimensional
subspace of $V$ such that $\varphi W \subseteq W$ and $\varphi^n V
\subseteq W$, for some $n$, then $\tr_V(\varphi) =
\tr_W(\varphi)$. And, since $\varphi$ is finite potent, we may
take $W = \varphi^n V$.\end{rem}

    This trace is the main tool used by J. Tate in his elegant
definition of the residue offered in \cite{Ta}.

    An open problem has been to determine whether this trace satisfies
the linearity property. In fact, in the article mentioned J. Tate
wrote: ``I doubt whether the rule
$$\tr_V \theta_1 + \tr_V \theta_2  = \tr_V (\theta_1 +
\theta_2)$$\noindent holds in general, i.e., whenever all three
endomorphisms $\theta_1$, $\theta_2$ and $\theta_1 + \theta_2$ are
finite-potent, although I do not know a counter-example. (If a
counter-example exists at all, then there will be one with
$\theta_1$ and $\theta_2$ nilpotent, because every finite-potent
endomorphism is the sum of a nilpotent one and one with finite
range.)''

    Lately, the second author of this note has studied this problem,
solving the question of the linearity with a negative answer for
the case of the sum of three finite potent endomorphisms.

    Indeed, in \cite{Pa} F. Pablos Romo offered three nilpotent
endomorphisms $\theta_1$, $\theta_2$ and $\theta_3$ of an
infinite-dimensional ${\mathbb Z}/2{\mathbb Z}$-vector space $V$,
such that $\theta_1 + \theta_2 + \theta_3$ is a finite-potent
endomorphism of $V$ and
$$\tr_V (\theta_1 + \theta_2 + \theta_3) = 1\, .$$

    Accordingly, in general Tate's trace does not satisfy the
linearity property.

    Moreover, M. Argerami, F. Szechtman and R. Tifenbach have
given an alternative characterization of finite potent
endomorphisms that can be used to reduce the question to a
special case. They have
shown in \cite{AST} that an endomorphism $\varphi$ is finite
potent if and only if $V$ admits a $\varphi$-invariant
decomposition $V = W_\varphi \oplus U_\varphi$ such that
$\varphi_{\vert_{U_\varphi}}$ is nilpotent, $W_\varphi$ is finite
dimensional, and $\varphi_{\vert_{W_\varphi}} \colon W_\varphi
\longrightarrow W_\varphi$ is an isomorphism. This decomposition
is unique and one has that $\tr_V(\varphi) =
\tr_W(\varphi_{\vert_{W_\varphi}})$.

    From this characterization, the question of the linearity of
Tate's trace can be reduced to a particular case of finite potent
linear operators in the vector space $V = k<s,t>/I$, where $I$ is
a left ideal of $k<s,t>$ such that:
\begin{itemize}
\item $I$ contains the two-sided ideal generated by $s^n, t^m,
f(s+t)$ for some $n,m \geq 2$ and $f(x) \in k[x]$, where $f(x) =
x^l g(x)$, $l\geq 2$ and $g(x) \ne 1$ is monic and relatively
prime to $x$.

\item If $J$ is the right ideal of $k<s,t>$ generated by
$(s+t)^l$, then $V_1 = (J + I)/I$ is finite dimensional and
non-trivial.

\item $V$ is infinite dimensional.
\end{itemize}

    The aim of this note is to give a negative answer to the question of the linearity of
Tate's trace for the sum of two endomorphisms by offering a counter-example of this property. For this counter-example
we consider a nilpotent endomorphism over a space V of countable dimension with diagonal entries 1,0,0,0,.... -exhibited in \cite{AST}-.

\section{Counter-example of the linearity property of the Tate's trace}

    Let $V$ be a vector space of countable dimension over an arbitrary ground field $k$.
Let $\{e_1,e_2,e_3,\dots\}$ be a basis of $V$ indexed by the natural numbers.

Let us consider the linear operator $\theta_1$ of $V$ defined in \cite{AST} by:

$$\theta_1 (e_i) = \left \{ {\begin{aligned} 0 \qquad &\text{ if i is odd } \\ e_{i-1} \,\ &\text{ if i is even } \end{aligned}} \right . \, .$$

If we denote by $\{v_j\}_{j\in {\mathbb N}}$ the basis constructed according to the following scheme:

$$v_1 = e_2, v_2 = e_2 + e_4, v_3= e_1 + e_4, v_{2i} = e_{2i} + e_{2i + 2}, v_{2i+1} = e_{2i-1} + e_{2i+2} \text{ for all } i\geq 1\, ,$$\noindent it is easy
to check that the matrix associated with $\theta_1$ in the basis $\{v_j\}_{j\in {\mathbb N}}$ is:

$$\theta_1 \equiv \begin{pmatrix} 1 & 0 & -1 & 0 & 1 & 0 & -1 & 0 & 1 &  \cdots \\
-1 & 0 & 1 & 0 & -1 & 0 & 1 & 0 & -1 &   \cdots \\ 1 & 1 & 0 & 0 & 0 & 0 & 0 & 0 & 0 &  \cdots \\ 0 & -1 & -1 & 0 & 1 & 0 & -1 & 0 & 1 &  \cdots \\
0 & 1 & 1 & 1 & 0 & 0 & 0 & 0 & 0 &  \cdots \\ 0 & 0 & 0 & -1 & -1 & 0 & 1 & 0 & -1 &  \cdots \\
0 & 0 & 0 & 1 & 1 & 1 & 0 & 0 & 0 &  \cdots \\ 0 & 0 & 0 & 0 & 0 & -1 & -1 & 0 & 1 &  \cdots \\
0 & 0 & 0 & 0 & 0 & 1 & 1 & 1 & 0 &  \cdots \\ \vdots & \vdots & \vdots & \vdots & \vdots & \vdots & \vdots & \vdots & \vdots &  \ddots  \end{pmatrix}\, .$$

Note that the explicit expression of $\theta_1$ in that basis is:

$$\theta_1 (v_i) = \left \{ {\begin{aligned} v_1 - v_2 + v_3 \qquad \qquad \qquad \qquad \qquad \qquad &\text{ if i = 1 } \\  v_3 - v_4 + v_5 \qquad \qquad \qquad \qquad \qquad \qquad &\text{ if i = 2 } \\
-v_1 + v_2 - v_4 + v_5  \qquad \qquad \qquad \qquad \qquad &\text{ if i = 3 } \\ v_{2k + 1} - v_{2k + 2} + v_{2k + 3} \qquad \qquad \qquad \qquad \qquad &\text{ if i = 2k for } k\geq 2
\\ (-1)^kv_1 + (-1)^{k+1}v_2 + [\sum_{\underset {j\geq 2} {2j = 4}}^{i-1} (-1)^{j+k}v_{2j}] - v_{i+1} + v_{i+2} \quad  &\text{ if i = 2k + 1 for } k\geq 2
 \end{aligned}} \right .$$

If $\{v_j\}_{j\in {\mathbb N}}$ is again the basis of $V$ described above, let us now consider the linear operator $\theta_2$ of $V$ defined by:

$$\theta_2 (v_i) = \left \{ {\begin{aligned} v_1 - v_3 \qquad \qquad \qquad &\text{ if i = 1 } \\ v_1 - v_2 - v_3 + v_4 - v_5 \qquad &\text{ if i = 2 } \\
v_1 - v_2 + v_4 \qquad \qquad &\text{ if i = 3 } \\ v_{4k + 2} \qquad \qquad \qquad &\text{ if i = 4k for } k\geq 1
\\ -v_1 + v_2 + \sum_{\underset {j\geq 2} {2j = 4}}^{i+1} (-1)^{j-1}v_{2j} \qquad &\text{ if i = 4k + 1 for } k\geq 1
\\ 0 \qquad \qquad \qquad \qquad &\text{ if i = 4k + 2 for } k\geq 1
\\ v_1 - v_2 + [\sum_{\underset {j\geq 2} {2j = 4}}^{i+1} (-1)^{j}v_{2j}] - v_{i+2} \qquad  &\text{ if i = 4k + 3 for } k\geq 1 \end{aligned}} \right .$$

A computation shows that $\theta_2$ is nilpotent of order 6 and its matrix with respect to the basis $\{v_j\}_{j\in {\mathbb N}}$ is:

$$\theta_2 \equiv \begin{pmatrix} 1 & 1 & 1 & 0 & -1 & 0 & 1 & 0 & -1 &  \cdots \\
0 & -1 & -1 & 0 & 1 & 0 & -1 & 0 & 1 &   \cdots \\ -1 & -1 & 0 & 0 & 0 & 0 & 0 & 0 & 0 &  \cdots \\ 0 & 1 & 1 & 0 & -1 & 0 & 1 & 0 & -1 &  \cdots \\
0 & -1 & 0 & 0 & 0 & 0 & 0 & 0 & 0 &  \cdots \\ 0 & 0 & 0 & 1 & 1 & 0 & -1 & 0 & 1 &  \cdots \\
0 & 0 & 0 & 0 & 0 & 0 & 0 & 0 & 0 &  \cdots \\ 0 & 0 & 0 & 0 & 0 & 0 & 1 & 0 & -1 &  \cdots \\
0 & 0 & 0 & 0 & 0 & 0 & -1 & 0 & 0 &  \cdots \\ \vdots & \vdots & \vdots & \vdots & \vdots & \vdots & \vdots & \vdots & \vdots &  \ddots  \end{pmatrix}\, .$$

Thus, if $\varphi = \theta_1 + \theta_2$, one has that $\varphi$ is a finite potent endomorphism of $V$, and the $\varphi$-invariant decomposition of $V$ referred to above is $V = W_\varphi \oplus U_\varphi$, with $W_\varphi = <v_1,v_2>$ and $U_\varphi = <v_r>_{r\geq 3}$.

 Regarding the basis $\{v_1,v_2\}$, it is clear that the isomorphism $\varphi_{\vert_{ W_\varphi}}$ is $$\varphi_{\vert_{ W_\varphi}} \equiv \begin{pmatrix} 2 & 1 \\ -1 & -1 \end{pmatrix}\, .$$

Accordingly, the explicit expression of the linear operator $\varphi_{\vert_{ U_\varphi}}$ in the basis $\{v_r\}_{r\geq 3}$ is:
$$\varphi_{\vert_{ U_\varphi}} (v_r) = \left \{ {\begin{aligned} v_5 \qquad \qquad \qquad &\text{ if r = 3 } \\  v_{4k + 1} + v_{4k + 3} \qquad \qquad &\text{ if r = 4k  for } k\geq 1
\\ v_{4k + 3} \qquad \qquad \qquad &\text{ if r = 4k + 1 for } k\geq 1
\\ v_{4k + 3} - v_{4k + 4} + v_{4k + 5}  \qquad  &\text{ if r = 4k + 2 for } k\geq 1
\\ 0 \qquad \qquad \qquad \qquad &\text{ if r = 4k + 3 for } k\geq 1 \end{aligned}} \right .$$
\noindent which is a nilpotent endomorphism of order 4, whose matrix in this basis is:

$$\theta_1 + \theta_2 \equiv \begin{pmatrix} 0 & 0 & 0 & 0 & 0 & 0 & 0 & 0 & 0 &  \cdots \\
0 & 0 & 0 & 0 & 0 & 0 & 0 & 0 & 0 &   \cdots \\ 1 & 1 & 0 & 0 & 0 & 0 & 0 & 0 & 0 &  \cdots \\ 0 & 0 & 0 & 0 & 0 & 0 & 0 & 0 & 0 &  \cdots \\
0 & 1 & 1 & 1 & 0 & 0 & 0 & 0 & 0 &  \cdots \\ 0 & 0 & 0 & -1 & 0 & 0 & 0 & 0 & 0 &  \cdots \\
0 & 0 & 0 & 1 & 0 & 1 & 0 & 0 & 0 &  \cdots \\ 0 & 0 & 0 & 0 & 0 & 0 & 0 & 0 & 0 &  \cdots \\
0 & 0 & 0 & 0 & 0 & 1 & 1 & 1 & 0 &  \cdots \\ \vdots & \vdots & \vdots & \vdots & \vdots & \vdots & \vdots & \vdots & \vdots &  \ddots  \end{pmatrix}\, .$$

Hence, $\theta_1$ and $\theta_2$ are nilpotent endomorphisms of $V$, and $\theta_1 + \theta_2$ is a finite potent endomorphism with $\tr_V (\theta_1 + \theta_2) = 1$. Thus, we obtain a counter-example of the linearity property of Tate's trace for finite
potent endomorphisms and we solve the above referred problem proposed by J. Tate in \cite{Ta}.

\bigskip {\centerline {\bf ACKNOWLEDGMENT}} The first author wishes to thank the ''Instituto Universitario de F\'isica Fundamental y Matem\'aticas (IUFFyM)" of the University of Salamanca for a Research Fellowship in the summer of 2011 that allowed her to explore this problem and begin an approach to its solution.

\end{document}